\documentclass{elsart}

\usepackage{amssymb}
\newcommand{\R}{{\Bbb R}}

\newcommand{\N}{{\Bbb N}}
\newcommand{\C}{{\Bbb C}}
\begin{document}
\begin{frontmatter}
\title{Non-monotone travelling waves
in a single species reaction-diffusion equation with delay}
\author[a]{Teresa Faria\corauthref{AA}}
\author[b]{and Sergei Trofimchuk}
\corauth[AA]{Corresponding author.}
\address[a]{Departamento de Matem\'atica, Faculdade de
Ci\^encias/CMAF, Universidade Lisboa, 1749-016, Lisboa, Portugal \\
{\rm E-mail: tfaria@ptmat.fc.ul.pt} }
\address[b]{Instituto de Matem\'atica y Fisica, Universidad de Talca, Casilla 747,
Talca, Chile \\ {\rm E-mail: trofimch@inst-mat.utalca.cl} }
\begin{abstract}
\noindent We prove the existence of a continuous family of
positive and generally non-monotone travelling fronts for delayed
reaction-diffusion equations $u_t(t,x) = \Delta u(t,x)- u(t,x) +
g(u(t-h,x)) \ (*)$, when  $g \in C^2(\R_+,\R_+)$ has exactly two
fixed points: $x_1= 0$ and $x_2= K
>0$. Recently, non-monotonic waves were observed
in numerical simulations by various authors. Here, for a wide
range of parameters, we explain why such waves appear naturally as
the delay $h$ increases. For the case of $g$ with  negative
Schwarzian, our conditions are rather optimal; we observe that the
well known Mackey-Glass type equations with diffusion fall within
this subclass of $(*)$. As an example, we consider the diffusive
Nicholson's blowflies equation.
\end{abstract}
\end{frontmatter}
\newpage
\section{Introduction}
\label{intro} In this paper, we study the existence of positive
non-monotone travelling waves for a family of delayed
reaction-diffusion equations which includes, as a particular case,
the diffusive Nicholson's blowflies equation
\begin{equation}\label{1s}
N_t(t,x)= d\Delta N(t,x) -\delta N(t,x) + p N(t-h,x)
e^{-bN(t-h,x)},
\end{equation}

\vspace{-5mm}

$t\in \R,\  x\in \R^m$. This problem was suggested   in
\cite{gou,gouk,gouss,lw,swz}, where numerical simulations
indicated a loss of monotonicity of the wave profile caused by the
delay. Eq. (\ref{1s}) was introduced in \cite{sy} and it
generalizes the famous Nicholson's blowflies equation
\begin{equation}\label{1ss}
N'(t)=  -\delta N(t) + p N(t-h) e^{-bN(t-h)},
\end{equation}

\vspace{-5mm}

intensively studied for the last decade (e.g. see our list of
references). After a linear rescaling of both variables $N$ and
$t$,  we can assume  that $\delta = b = 1$. Eq. (\ref{1s}) takes
into account spatial distribution of the species, and the
mentioned problems reflect the  interest in understanding
 the spatial spread of the growing population \cite{gn}. Relevant
biological discussion can be found in
\cite{gous,gou,gour,gouss,lw,swz}, where various modifications of
(\ref{1s}) were proposed and studied. Here, however, we will
concentrate mainly on the mathematical aspects of the dynamics in
(\ref{1s}). For the sake of simplicity, we will consider the case
of a single discrete delay, but extensions for more general
functionals (which additionally can take into account non local
space effects) are possible (cf. \cite{fhw,ltt1,ltt}).  Since the
biological interpretation of $N$ is the size of an adult
population, we will consider {\it only} non-negative solutions for
(\ref{1s}) and for other population models. Actually, our approach
allows us to study a more general
family of scalar reaction-diffusion equations 
\begin{equation}\label{17}
u_t(t,x) = d\Delta u(t,x)  - u(t,x) + g(u(t-h,x)), \ u(t,x) \geq
0,\ x \in \R^m,
\end{equation}

\vspace{-5mm}

related to the Mackey-Glass type delay differential equations,
\begin{equation}\label{17a}
u'(t) =   - u(t) + g(u(t-h)), \ u \geq 0,
\end{equation}

\vspace{-5mm}

with exactly two non-negative equilibria $u_1(t) \equiv 0, \
u_2(t) \equiv K >0$ (so that $g(K) = K, \ g(0) =0$). In
particular, with $g(u)= pu/(1+u^n)$ in (\ref{17a}), we obtain the
equation proposed in 1977 by Mackey and Glass, to model
hematopoiesis (blood cell production). The non-linearity $g$ is
called {\it the birth function} and thus it is non-negative, and
generally non-monotone and bounded. Due to these properties of $g$
and the simple form of dependence on the delay in (\ref{17}), the
Cauchy problem
\begin{equation} \label{CP} u(s,x) = \zeta (s,x), \ s \in [-h,0],
\ x \in \R^m,
\end{equation}
for Eq. (\ref{17a}) has a unique eventually positive global
solution for every $\zeta\not=0$ taken from an appropriately
chosen functional space (e.g. see \cite{smZH}).

Recently, the existence of  travelling fronts connecting the
trivial and positive steady states in $(\ref{1s})$ was studied in
\cite{sz} and \cite{fhw} (see also \cite{gou,ma} for other methods
which eventually can be applied to analyze this problem). In
\cite{sz}, the authors use a monotone iteration procedure coupled
with the method of upper and lower solutions. This approach
(proposed in \cite{wz}) works well if $1<p/\delta \leq e,$ since
in this case the function $g$ is increasing on $[0,1/b]\supset
[0,K]$, thus $\phi \mapsto -\delta \phi (0)+p g(\phi (-h))$
satisfies the quasimonotonicity condition in \cite{wz}. This
allows one to establish the existence of {\it monotone} wave front
solutions $N(t,x) = \phi(ct + \nu\cdot x, c)$ for every $p/\delta
\in (1, e]$  and  $c
> 2 \sqrt{p-\delta}$ (cf. \cite{mei,sz}). Moreover, as it was
proved in \cite{mei}, every solution of $(\ref{1s}),(\ref{CP})$
with $p/\delta \in (1, e]$  converges exponentially to some
travelling wave provided that $\zeta$ is sufficiently close (in a
weighted $L^2$ norm) to this wave at the very beginning of the
propagation. For the case $p/\delta
>e$, clearly $g$ is not monotone on $[0,K]$, and  Wu and Zou's
method \cite{wz} is no longer applicable.  In \cite{fhw}, the
Lyapunov-Schmidt reduction was used to study systems of delayed
reaction-diffusion equations with non-local response. We observe
that Eqns. (\ref{1s}) and (\ref{17})  fit into the framework
developed in \cite{fhw}. This approach requires a detailed
analysis of an associated Fredholm operator and the existence of
heteroclinic solutions of (\ref{17a}) (in \cite{fhw}, the latter
was established with the use of the monotone semiflows approach
developed by H. Smith and H. Thieme \cite{SmTh}). As a result, it
was proved in \cite{fhw} that, even when $p/\delta > e$,
(\ref{1s}) possesses a family of travelling waves if $\delta h \in
(0, r^*)$ for some $r^* < 1$ (which is given explicitly). The
rather restrictive condition $\delta h < r^* < 1$ from \cite{fhw}
was considerably weakened in \cite{ft} by invoking a Schauder's
fixed-point argument to find heteroclinic solutions of
(\ref{17a}). Unfortunately, the main results of \cite{fhw,ft} do
not answer the question about the existence (and shape) of {\it
positive} travelling fronts of (\ref{1s}) or (\ref{17}). We recall
here that only non-negative solutions to (\ref{17}) are
biologically meaningful.

In this paper, inspired by \cite{fhw,sz,wz},  for a broad family
of nonlinearities $g$ (which includes Eq. (\ref{1s}) with $\delta
=1$),  we prove that Eq. (\ref{17})  has a continuous family of
positive travelling wave fronts $u(t,x) = \phi(ct + \nu\cdot x,
c)$, indexed by the speed parameter $c
> 0$,  provided that 
\begin{equation} \label{gsc}
e^{- h} > - \Gamma \ln \frac{{\Gamma}^{2} - \Gamma}{{\Gamma}^{2} +
1}, \quad \Gamma \stackrel{def}= g'(K),
\end{equation}

\vspace{-5mm} and $c$ is sufficiently large: $c >
c_*(h,g'(0),g'(K))$. Furthermore, we show that these fronts
generally are not monotone: in fact, they can oscillate infinitely
about the positive steady state. On the other hand, for large
negative values of $s$, the wave profile $\phi(s, c)$ is
asymptotically equivalent to an increasing exponential function.
Condition (\ref{gsc}) assures the global attractivity of the
positive equilibrium of  (\ref{17a}), which is required by our
approach. It should be noted that this condition is rather
satisfactory in the sense that (\ref{gsc})  determines a domain of
parameters approximating very well the maximal region of local
stability for the positive steady state in (\ref{17a}) or
(\ref{1ss}) (cf. \cite{ltt}).

Before announcing the main results of the present work, we state
our basic hypothesis:
\begin{description}
\item[{\rm \bf(H)}] Eq. (\ref{17a}) has exactly two steady states
$u_1(t) \equiv 0$ and $u_2(t)\equiv K>0$, the second equilibrium
being exponentially asymptotically stable and the first one being
hyperbolic.  Furthermore, $g \in C^1(\R_+, \R_+)$ and is $C^2$-
smooth in some vicinity of the equilibria, with $p:=g'(0)>1$. The
latter implies that the solution $u_1=0$ of (\ref{17a}) is
unstable for all $h\ge 0$.
\end{description}
In the sequel,  $\lambda_1(c)$  denotes the minimal positive root
of the characteristic equation $(z/c)^2-z-1+p\exp(-z h)=0$ for
sufficiently large $c$, and $\lambda$  the unique positive root of
the equation $-z-1+p\exp(-z h)=0$, where $p > 1$. As shown later,
$\lim_{c \to \infty} \lambda_1(c) = \lambda$. Now we are ready to
state our main result:
\begin{thm} \label{main} Assume {\rm \bf(H)}. If the positive equilibrium $K$
of Eq. (\ref{17a}) is  globally attracting,  then there is $c_*
> 0$ such that, for each  $\nu \in \R^m, \ \|\nu\| =1$, equation (\ref{17}) has a continuous family of positive
travelling waves $u(t,x) = \phi(ct + \nu\cdot x, c), \ c > c_*$.
Furthermore, for some $s_0 = s_0(c) \in \R$, we have $\phi(s-s_0,
c) = \exp (\lambda_1(c) s) + O(\exp(2\lambda s))$ as $s \to -
\infty$, so that $\phi'(s-s_0,c) = \lambda_1(c) \exp (\lambda_1(c)
s) + O(\exp(2\lambda s))> 0$ on some semi-axis $(-\infty, z]$.
Finally, if $\ g'(K) he^{h+1} <-1$ then the travelling profile
$\phi(t)$ oscillates about $K$ on every interval $[z, +\infty)$.
\end{thm}
In order to apply Theorem \ref{main}, one needs to find sufficient
conditions to ensure the global attractivity of the positive
equilibrium of (\ref{17a}). Some results in this direction were
found in \cite{ltt1,ltt} for a family of nonlinearities having
negative Schwarz derivative (or, more generally, satisfying a
generalized Yorke condition \cite{flot,ltt1,ltt}). In particular,
\cite[Corollary 2.3]{ltt1} implies the following useful version of
Theorem \ref{main}: \noindent
\begin{cor} \label{nasa} Assume {\rm \bf(H)} and (\ref{gsc}).
In addition,  suppose  that $g \in C^3(\R_+, \R_+)$ has only one
critical point $x_M$ (maximum) and  that the Schwarz derivative
$(Sg)(x)=g'''(x)(g'(x))^{-1}-(3/2)
\left(g''(x)(g'(x))^{-1}\right)^2 $ is negative for all $x >0$, $x
\not=x_M$. Then all conclusions of Theorem 1 hold true.
\end{cor}
Notice that Corollary \ref{nasa} applies to both the Nicholson's
blowflies equation and the Mackey-Glass equation with non-monotone
nonlinearity, see \cite{ltt1}.

To prove our main results, we need a detailed analysis of
heteroclinic solutions of (\ref{17a}). This study is presented in
Section 2, and is crucial for the selection of an appropriate
functional space where a Lyapunov-Schmidt reduction is realized.
The existence of positive travelling waves is proven in the third
section. The main result of Section 3 is given in Theorem
\ref{mr}, which is essentially Theorem \ref{main} without its
non-monotonicity statement. Finally, in the last short section, we
show that these waves have non-monotonic profiles when the delay
is over some critical value.
 \vspace{-5mm}
\section{Heteroclinic solutions of scalar delay differential equations}
\vspace{-2mm}
 In this section, we study the existence and
properties of heteroclinic solutions to the scalar functional
equation
\begin{equation} \label{1}
x'(t)= - x(t) + f(x_t), \ x \geq 0,
\end{equation}
\vspace{-4mm} where $f: C([-h,0], \R_+) \to \R_+$ is a continuous
functional which takes closed bounded sets into bounded subsets of
$\R_+$. Here $C([-h,0], \R_+)$ is the metric space equipped with
the norm $|\phi| = \max_{s \in [-h,0]}|\phi(s)|$.  Throughout this
paper, we suppose that Eq. (\ref{1}) has exactly two steady states
$x_1(t) \equiv 0$ and $x_2(t)\equiv K$, the second equilibrium
being asymptotically stable and globally attractive. Thus, if
(\ref{1}) has a heteroclinic solution $\psi(t)$, it must satisfy
$\psi(-\infty )=0, \ \psi(+\infty) = K$.

We start by proving a general existence result which is valid for
the abstract setting of dynamical systems. Let $S^t:X \to X$ be a
continuous semidynamical system defined in a complete metric space
$(X,d)$. First, we mention the following fact (see e.g.  \cite[p.
36]{HMO}):
\begin{lem}
Suppose that $\varphi:\R\to X, \ \varphi(0)  =x$ is a complete
orbit of $S^t$.  If the closure $\overline{\{\varphi(s), \ s \leq
p \}}$ is compact for some $p\in \R$, then the $\alpha-$limit set
$\alpha(\varphi) =\cap_{q \leq 0}
 \overline{\{\varphi(s), \ s \leq q \}}$ of $\varphi$ is nonempty, compact and
invariant (this means that for every $z \in \alpha(\varphi)$ there
exists at least one full trajectory  $\psi$ with $\psi(\R)
\subseteq \alpha(\varphi), \ \psi(0)=z$).
\end{lem}
For every $A \subset X$ and $h > 0$, let $A(h) \subset A$ denote
the set of right endpoints of all orbit segments $S^{[0,h]}z =
\{S^uz: u\in [0,h]\}$ which are completely contained in $A$:
$$
A(h) = \{x \in A: x = S^hz \, {\rm and }\ S^{[0,h]}z \subset A,  \
{\rm for \ some \ } z \in A\}.
$$
Next statement shows clearly how to relate the global attractivity
property of the positive equilibrium of (\ref{17a}) to the problem
concerning the existence of travelling fronts for (\ref{17}):
\begin{lem} \label{123} Assume that $A(h)$ is either empty or pre-compact, for all
bounded sets $A$ and  some $h > 0$. Suppose that there exist two
disjoint compact invariant subsets $K_1, K_2$ of $X$ such that
$d(S^tx,K_2) \to 0$ as $t \to +\infty$ for every $x \in X\setminus
K_1$. If the set $F_{\varepsilon}= \{x: d(x,K_1) = \varepsilon \}$
is not empty for every sufficiently small $\varepsilon > 0$, then
there exists at least one complete orbit $\psi$ with $\alpha(\psi)
\subset K_1$ and $\omega(\psi) \subset K_2$.
\end{lem}
\begin{pf} Let $\rho = d(K_1,K_2)$ and, for every $n > 2/\rho$, take some
$x_n \in F_{1/n}$. Due to the compactness of $K_1$,  we can assume
that $x_n \to z$ for some $z \in K_1$. In consequence, if $t_n
>0$ is the minimal real number such that $d(S^{t_n}x_n, K_1) = \rho/2$,
then $\lim t_n = +\infty$. Set $w_n = S^{t_n}x_n$. Due to the
compactness condition imposed on $S^t$, we can suppose that $\lim
w_n = w$. Let now $\psi_n(u)= S^{u+t_n}x_n, \ u \geq -t_n$. We
have $S^a\psi_n(t) = \psi_n(a+t)$ for every $a \geq 0, t \geq
-t_n$. Since, for every integer $m > 0$ the sequence $\psi_n(-m)$
has a convergent subsequence (say, $\psi_{n_{j}}(-m) \to b$), we
can assume that $\psi_n(t)$ converges uniformly on $[-m,0]$ to
$\psi(t) = S^{t+m}b$. Moreover, we have that $\psi(0)=w$ and
$S^a\psi(t) = \psi(a+t)$ for all $a \geq 0, t \geq -m$. In this
way, taking $m= 1,2,3, \dots$, we can use $\psi_n(u)$ to construct
a continuous function $\psi: \R \to X$, such that $S^a\psi(t) =
\psi(a+t)$ for every $a \geq 0,\ t \in \R$. Such $\psi$ defines
the complete orbit we are looking for. Since $\psi(\R_{-})$ is a
subset of the bounded set $B = \{z: d(z, K_1) \leq \rho/2\}$, we
conclude that $\psi(\R_{-})$ is pre-compact. Furthermore, because
of $d(\psi(\R_{-}), K_2) \geq \rho/2$, we have
$d(\alpha(\psi),K_2) \geq \rho/2
>0$. This means that $\alpha(\psi) \subset K_1$ so that
$\lim_{t \to - \infty}d(\psi(t), K_1) =0$. $\qquad \Box$
\end{pf}

\vspace{-5mm}

A direct application of Lemma \ref{123} to  Eq. (\ref{1}) gives
the following
\begin{thm} Let $f: C([-h,0], \R_+) \to \R_+$ be a continuous functional
which takes closed bounded sets into bounded subsets of $\R_+$.
Assume further that every non-negative solution of (\ref{1})
admits a unique extension on the right semi-axis. If $f(0)=0,$ $
f(K)=K\ (K>0)$ and $x_2(t) \equiv K$ attracts every solution of
(\ref{1}) with non-negative and nontrivial initial function, then
there exists a positive complete solution $\psi(t)$ to (\ref{1})
such that $\psi(- \infty) = 0, \ \psi(+ \infty) = K$.
\end{thm}
With some additional conditions on $f$, we can say more about such
an orbit $\psi$:
\begin{lem} \label{4} Assume $f(\phi) = g(\phi(-h))$ for some $g \in C(\R_+,\R_+)$ (so that
$g(0)=0, \ g(K) = K$). Assume that $\liminf_{x \to 0+} g(x)/x>1$,
and let $p_1,p_2$ be such that $1<p_1 < \liminf_{x \to 0+}
g(x)/x\leq \limsup_{x \to 0+} g(x)/x<p_2$. Let $\lambda_i$
be the unique positive real root of the 
equation $z = - 1 + p_i\exp (-z h)$ ($i=1,2$), so we have $0<
\lambda_1<\lambda_2$.  Then for every heteroclinic solution
$\psi(t)$ of the equation
\begin{equation} \label{2}
x'(t)= - x(t) + g(x(t-h))
\end{equation}

\vspace{-5mm}

there exist $\tau = \tau(\psi) < 0,\ C_i = C_i(\psi)> 0$ such that
$$
C_1\exp(\lambda_2 t) \leq \psi (t) \leq C_2\exp(\lambda_1 t), \ t
\leq \tau.
$$
\end{lem}

\vspace{-7mm}

\begin{pf}
Choose $\delta >0$  sufficiently small such that  $p_1x \leq g(x)
\leq p_2x$ for all $x \in [0,\delta)$, and let $\tau$ be such that
$\psi(t) < \delta$ for all $t \leq \tau$. We claim that, for every
$s \leq \tau$
\begin{equation}\label{es} \hspace{-7mm}
\psi(t_m)=\min_{u \in [-h,0]} \psi(s+u) \geq
\frac{\exp(-h)}{p_2}\max_{u \in [-h,0]} \psi(s+u)
=\frac{\exp(-h)}{p_2}\psi(t_M).
\end{equation}
Indeed, if $t_m \geq t_M$, then, by the variation of constants
formula,
$$
\psi(t_m)  = \psi(t_M)\exp(t_M-t_m) +
\int_{t_M}^{t_m}\exp(-(t_m-u))g(\psi(u-h))du \geq
\psi(t_M)\exp(-h). $$ Finally, suppose that $t_M-h \leq t_m < t_M$
so that  $\psi'(t_M) \geq 0$. Then (\ref{es}) holds since
\begin{eqnarray*}\label{Y} \hspace{-7mm}
&  & \psi(t_M) \leq g(\psi(t_M-h)) \leq p_2\psi(t_M-h); \
\psi(t_m) = \psi(t_M-h)\exp(t_M-t_m-h) + \\ \hspace{-7mm} &
&\int_{t_M-h}^{t_m}\exp(-(t_m-u))g(\psi(u-h))du \geq
\psi(t_M-h)\exp(-h).
\end{eqnarray*}
Next, for every $s \leq \tau$ and $u \in [-h,0]$, we have that
$$
\psi(t_m)\exp(\lambda_1 u)\leq \psi(s+u) \leq
\psi(t_M)\exp(\lambda_2(u+h))\leq
p_2\psi(t_m)\exp(\lambda_2(u+h)+h).
$$
From the inequalities above and since additionally $p_1x \leq g(x)
\leq p_2x$ for all $x \in [0,\delta)$,  then for $s+h\leq \tau$
and $u\in [-h,0]$ we have
\begin{eqnarray*}
& & \psi (s+h+u)= \psi (s)e^{-(h+u)}+e^{-(s+h+u)}\int_s^{s+h+u} e^{\sigma} g(\psi (\sigma -h))d\sigma \nonumber \\
& \leq & \psi (s)e^{-(h+u)}+e^{-(s+h+u)} p_2 \psi (t_M) \int_s^{s+h+u}\e^{\sigma+\lambda_2(\sigma -s)} d\sigma \nonumber \\
& = & \psi (s)e^{-(h+u)}+ e^{\lambda_2 h}  \psi (t_M)
[e^{\lambda_2(h+u)}-e^{-(h+u)}]  \leq  \psi (t_M) e^{\lambda_2
(u+2h)}, \nonumber
\end{eqnarray*}
and
\begin{eqnarray}
\psi (s+h+u)& \geq & \psi (s)e^{-(h+u)}
+e^{-(s+h+u)} p_1 \psi (t_m) \int_s^{s+h+u} e^{\sigma+\lambda_1(\sigma -h-s)} d\sigma \nonumber \\
& = & \psi (s)e^{-(h+u)}+  \psi (t_m)
[e^{\lambda_1(h+u)}-e^{-(h+u)}]
 \geq  \psi (t_m) e^{\lambda_1 (u+h)} .\nonumber
\end{eqnarray}
By repeating the above procedure over intervals of length $h$, the
step by step method  implies that, for all $-h \leq u \leq \tau
-s$,
$$
\psi(t_m)\exp(\lambda_1 u)\leq \psi(s+u) \leq
p_2\psi(t_m)\exp(\lambda_2(u+h)+h).
$$
In particular,
$$
\frac{\psi(s)}{p_2}\exp(\lambda_1 (\tau-s)-h)\leq
\psi(t_m)\exp(\lambda_1 (\tau-s))\leq \psi(\tau) \leq
p_2\psi(t_m)\exp(\lambda_2(\tau-s+h)+h).
$$
Thus, for every $s \leq \tau$,
$$
p_2^{-1}\psi(\tau)\exp(\lambda_2 (-\tau - h+ s)-h) \leq \psi(t_m)
\leq \psi(s)\leq p_2\psi(\tau)\exp(\lambda_1 (-\tau+ s)+h). \quad
\Box
$$

\end{pf}

In what follows, we shall assume that  $g \in C^1(\R_+, \R_+)$,
$g'(0+) = p >1$, and use several times the following simple
assertion.
\begin{lem}\label{char3}
Suppose that $p> 1$ and  $h >0$. Then the characteristic equation
\begin{equation}\label{char}
z = - 1 + p\exp (-z h)
\end{equation}
has only one real root $0 <  \lambda < p-1$. Moreover, all roots
$\lambda, \lambda_j, \ j =2,3,\dots $ of (\ref{char}) are simple
and we can enumerate them in such a way that $$\dots \leq \Re
\lambda_3 = \Re \lambda_2 < \lambda.$$
\end{lem}
\begin{pf} The last inequality follows from
$\Re\lambda_j < -1 + p \exp (- h \Re\lambda_j), \ j >1$. $\Box$
\end{pf}
\begin{lem} \label{psi} Suppose that $g'(0+) = p >1$ and that $\psi$
is a heteroclinic solution  to (\ref{2}). Let $\lambda$ be the
positive root of (\ref{char}). If there exists $g''(0+) \in \R$,
then, for each $\delta > 0$ and some $t_0 \in \R$, we have that
$\psi(t-t_0) = \exp (\lambda t) + O(\exp((2\lambda- \delta)t))$ at
$t \to - \infty$, so that $\psi'(t-t_0) = \lambda \exp (\lambda t)
+ O(\exp((2\lambda- \delta)t))> 0$ on some semi-axis $(-\infty,
T]$. Moreover, if there exists  $g'':[0,\epsilon )\to \R$ and is
bounded for some $\epsilon > 0$,  then $\psi(t)$ is unique up to a
shift in $t$.
\end{lem}

\vspace{-5mm}

\begin{pf}
Since $g''(0+)$ is finite, $g(x) = g'(0+)x + O(x^2)= px +  O(x^2)$
as $x \to 0$. From Lemma 6, given $\delta >0$ small, for
$p_1=p-\delta /2$ we have $\psi (t)=O(\exp
 (\lambda_\delta t))$ at $-\infty$, where $\lambda_\delta$ is the unique
 positive root of $z=-1+(p-\delta/2)\exp(-zh)$. It is easy to see  that
 $\lambda_\delta >\lambda-\delta/2$. In fact, let
 $W(z):=z+1-(p-\delta/2)\exp
 (-zh)$. For $z\in\R$, we have $W(z)<0$ if and only if $z<\lambda_\delta$. On the other
 hand, $W(\lambda -\delta/2)=p\exp (-\lambda h)[1-\exp(\delta
 h/2)]+\delta/2[\exp(-(\lambda-\delta/2)h)-1]<0.$
Hence,
\begin{equation}\label{ps}
\psi'(t) = - \psi(t) + p \psi (t-h) + O(\psi^2(t-h)),
\end{equation}
where $O(\psi^2(t-h))= O(\exp((2\lambda- \delta)t))$ as $t \to
-\infty$. Now, consider the linear inhomogeneous delay
differential equation
\begin{equation}\label{lps}
x'(t) = - x(t) + px(t-h) + O(\exp((2\lambda- \delta)t)).
\end{equation}
The change of variables $x(t) = y(t)\exp((2\lambda- \delta)t)$
transforms it into
\begin{equation}\label{llps}
y'(t) = -(1+2\lambda- \delta) y(t) + p\exp(-(2\lambda-
\delta)h)y(t-h) + O(1).
\end{equation}
The spectra $\sigma(y), \ \sigma(x)$ of the linear parts of
Eqns. (\ref{llps}) and  (\ref{lps}) are related by $\sigma(y)=
\sigma(x) -2\lambda + \delta$, therefore the linear part of
(\ref{llps}) has not pure imaginary eigenvalues for all
sufficiently small $\delta > 0$ (equivalently, the linearization
of Eq. (\ref{llps}) about zero  is hyperbolic).  In this case,
(\ref{llps}) has a bounded solution $y_b(t) = O(1)$ at $t = -
\infty$ (e.g. see \cite[Lemma 3.2, p. 246]{dglw} or \cite[Section
10.1]{hl}. Note that  Eqns. (\ref{lps}) and (\ref{llps}) are not
autonomous. Nevertheless, the results for autonomous equations
near hyperbolic equilibria in \cite{dglw,hl} are valid in this
setting, since the linearized equation near zero for Eq.
(\ref{llps}) has an exponential dichotomy, cf. \cite[p. 312]{hl}).
Thus Eq. (\ref{lps}) has a solution $x_b(t) =
y_b(t)\exp((2\lambda- \delta)t) = O(\exp((2\lambda- \delta)t)). $
In consequence, $z(t) = \psi(t) - x_b(t)$ solves the linear
homogenous equation $ x'(t) = - x(t) + px(t-h) $ and is bounded at
$t \to - \infty$. This is possible if and only if
$$
z(t) = C \exp (\lambda t) + \sum_{j=1}^NC_j \exp (\lambda_j t),
$$
where $\lambda > 0, \ \lambda_j \in \C, \ j = 1, \dots, N$ is a
finite set of roots having non-negative real parts of the
characteristic equation (\ref{char}). Notice that $C \in \R, C_j
\in \C$ and $\lambda > \Re \lambda_j$ (see Lemma \ref{char3}). In
this way
$$
\psi (t) = C \exp (\lambda t) + \sum_{j=1}^NC_j \exp (\lambda_j t)
+ O(\exp((2\lambda- \delta)t)).
$$
On the other hand, from Lemma \ref{4} we know that $\psi (t) =
O(\exp((\lambda- \delta/2)t))$. Since $\lambda > \Re \lambda_j$,
this implies immediately that all $C_j =0, \ C > 0$ and that $
\psi (t) = C \exp (\lambda t) + O(\exp((2\lambda- \delta)t)).$ By
(\ref{ps}),
$$
\psi'(t) = C\lambda \exp (\lambda t) + O(\exp((2\lambda-
\delta)t)) > 0.
$$
Observe also that $ \mu(t) = \psi (t - \lambda^{-1}\ln C) = \exp
(\lambda t) + O(\exp((2\lambda- \delta)t))$ defines another
heteroclinic solution of (\ref{1}).

Finally, suppose that $\mu(t), \nu(t)$ are two heteroclinic
solutions to (\ref{1}) such that
$$
\mu(t) = \exp (\lambda t) + O(\exp((2\lambda- \delta)t)), \ \nu(t)
= \exp (\lambda t) + O(\exp((2\lambda- \delta)t)).
$$
Applying the Lagrange mean value theorem twice, we  get
$g(x)-g(y)=  p(x-y) + (x-y)O(x+y)$ for $x,y$ close to $0$. Since
$\sigma(t) = \mu(t)-\nu(t) = O(\exp((2\lambda- \delta)t))$ we
obtain that
$$
g(\mu(t-h))- g(\nu(t-h)) = \sigma(t-h)(p + O(\exp (\lambda t))) =
p\sigma(t-h) + O(\exp ((3\lambda - \delta) t)).
$$
Therefore $\sigma(t)$ satisfies
\begin{equation}\label{alps}
x'(t) = - x(t) + p\, x(t-h) + O(\exp((3\lambda- \delta)t)),
\end{equation}

\vspace{-5mm}

from which, applying the same procedure as above, we deduce that
$\sigma(t) = \mu(t)-\nu(t) = O(\exp((3\lambda- \delta)t))$. In
this way, we can show that $\sigma(t) = O(\exp((k\lambda-
\delta)t))$ for every integer $k \geq 2$. This leads us to the
conclusion that $\sigma$ has superexponential decay at $t = -
\infty$ (equivalently, $\sigma$ is a {\it small solution} at $t =
- \infty$, see \cite{dglw}). We will finalize our proof showing
that only the trivial solution of the linear asymptotically
autonomous homogeneous equation
\begin{equation}\label{se}
x'(t) = - x(t) + p(t)x(t-h), \quad p(-\infty) = p > 1,
\end{equation}

\vspace{-5mm}

can have superexponential decay at $t = - \infty$ (notice that
$\sigma(t)$ satisfies (\ref{se}) with $p(t) = p + O(\exp (\lambda
t))$). Indeed, if $x(t)
> 0$ on some semi-axis $(-\infty, z]$, then we can repeat the
arguments in the proof of Lemma \ref{4} to find an exponential
lower bound for $x(t)$, in contradiction to our assumption of
superexponential decay of $x(t)$. Consider now the case of $x(t)$
oscillatory on every semi-axis $(-\infty, z]$, and take $z_0$ such
that $p(t)< p+1$ for all $t \in (-\infty, z_0]$. Let  $t_1 \in
(-\infty, z_0)$ be a point of the global maximum of $|x(t)|$: we
can assume that $x(t_1) =M > 0, \ x'(t_1) \geq 0$. Then $x(t_1-h)
\geq M/(p+1)$, so that $|x(t_2)| =  \max_{t \leq t_1-h}|x(t)| \geq
(p+1)^{-1}\max_{t \leq t_1}|x(t)|$. Analogously, $\max_{t \leq
t_1-2h}|x(t)| \geq |x(t_3)| =  \max_{t \leq t_2-h}|x(t)| \geq
(p+1)^{-1}\max_{t \leq t_2}|x(t)|= (p+1)^{-1}\max_{t \leq
t_1-h}|x(t)| \geq (p+1)^{-2}\max_{t \leq t_1}|x(t)|$. Thus
$$
\max_{t \leq t_1-kh}|x(t)|  \geq (p+1)^{-k}\max_{t \leq t_1}|x(t)|
$$
so that $x(t)$ can not decay superexponentially  as $t \to -
\infty$. $\qquad \Box$
\end{pf}

Now, assume {\rm \bf(H)} and the global attractivity of $x_2=K$
for Eq. (\ref{2}), and then take $\lambda >0$ satisfying
(\ref{char}) and the unique (up to a shift in time) heteroclinic
solution $\psi$ described in  Lemma \ref{psi}. Let $\lambda_* \in
(0, \lambda)$ be sufficiently close to $\lambda$  and such that
the equation $y'(t) = -(1+\lambda_*) y(t) +
p\exp(-\lambda_*h)y(t-h)$ is hyperbolic. Note that this latter
equation is obtained by effecting the change of variables
$x(t)=\exp (\lambda_* t)y(t)$ to the linear equation
$x'(t)=-x(t)+px(t-h)$.  For a fixed $\mu
> 0$,   we will consider the seminorms
$\|x\|^+ = \sup_{\R_+}|x(s)|$, $\|x\|^-_{\mu} = \sup_{\R_-}e^{-\mu
s}|x(s)|$, $\|x\|_{\mu} = \max\{\|x\|^+,\|x\|^-_{\mu}\}$ and the
following Banach spaces:
\begin{eqnarray*}\label{BS}
&  & C_{\mu}(\R) = \{x \in C(\R,\R): \|x\|^-_{\mu} < \infty {\rm \
and \ } x(+\infty) \ {\rm exists \ and \ is\  finite} \},
 \\
&  & C_{\psi,\lambda_*}(\R) = \{x \in C_{\lambda_*}(\R):
\int_{-\infty}^0 x(s)\psi'(s)ds =0 \},
\end{eqnarray*}
equipped with the norms $\|x\|_{\mu}$ and $\|x\|_{\lambda_*}$
respectively (in order to simplify the notation, we shall often  write $\|x\|$
instead of $\|x\|_{\mu}$, etc). Notice that, due to Lemma
\ref{psi}, we have $\psi, \psi' \in C_{\lambda_*}(\R)\setminus
C_{\psi,\lambda_*}(\R)$. We shall also need the following integral
operator 
$$
\mathcal{N}: C_{\psi,\lambda_*}(\R) \to C_{\lambda_*}(\R); \quad
(\mathcal{N}x)(t) = \int^t_{-\infty}e^{-(t-s)}q(s)x(s-h)ds,
$$
where $q(s) = g'(\psi(s-h))$ with $q(-\infty) = g'(0+)= p > 1, \
q(+\infty) = g'(K)$. Observe that $\mathcal{N}$ is well defined,
since $(\mathcal{N}x)(+\infty) = g'(K)x(+\infty)$ and, for $t \leq
h$,
$$
|(\mathcal{N}x)(t)| =
\int^t_{-\infty}e^{-(t-s)}|q(s)|\|x\|^-_{\lambda_*}e^{\lambda_*(s-h)}ds
\leq \frac{\|x\|^-_{\lambda_*}\sup_{t \le
h}|q(t)|}{1+\lambda_*}e^{\lambda_*
(t-h)}.
$$
\begin{lem} \label{6} If  {\rm \bf(H)} is assumed, then $I - \mathcal{N}:
C_{\psi,\lambda_*}(\R) \to C_{\lambda_*}(\R)$ is an isomorphism of
Banach spaces.
\end{lem}
\begin{pf}
We first prove that ${\rm Ker}(I - \mathcal{N})=0$. Indeed, if $ y
\in {\rm Ker}(I - \mathcal{N})$ and $y\not= 0$, then
$$
\int^t_{-\infty}e^{-(t-s)}q(s)y(s-h)ds = y(t).
$$
Therefore $y$ is a bounded solution of the linear delay
differential equation
\begin{equation}\label{veq} y'(t) = - y(t) +
q(t)y(t-h).
\end{equation}
Since $g'(x) = p+ O(x)$ at $x=0$ and $\psi(t) = O(\exp(\lambda
t))$ at $t = -\infty$, we conclude that
$$
q(t) = g'(\psi(t-h)) = p + O(\exp(\lambda t)), \ t \to -\infty.
$$
Thus $y(t)$ can be viewed as a bounded solution of the
inhomogeneous equation
$$
x'(t) = - x(t) + px(t-h) + O(\exp(2\lambda t)), \ t \to -\infty.
$$
Since $y(t) = O(\exp(\lambda_* t))$ at $-\infty$, with $\lambda_* <\lambda$
close to $\lambda$, the procedure which has been used before to
prove the uniqueness of the heteroclinic $\psi(t)$ allows us to
conclude that $y(t) = C\exp(\lambda t)+ O(\exp(2\lambda_* t))$ and
that $\dim{\rm Ker}(I - \mathcal{N})=1$. On the other hand, we
know that $\psi'(t)\not\equiv 0$ satisfies (\ref{veq}). Thus we
must have $y(t) = c\psi'(t)\not\in C_{\psi,\lambda_*}(\R)$, $c\ne
0$ constant, a contradiction. Therefore $y(t) \equiv 0$ and ${\rm
Ker}(I - \mathcal{N})=0$.

We now  establish that $I - \mathcal{N}$ is an epimorphism. Take
some $d \in C_{\lambda_*}(\R)$ and consider the following integral
equation
$$
x(t) - \int^t_{-\infty}e^{-(t-s)}q(s)x(s-h)ds = d(t).
$$
If we set $z(t) = x(t)-d(t)$, this equation is transformed into
$$
z(t) - \int^t_{-\infty}e^{-(t-s)}q(s)(z(s-h)+d(s-h))ds = 0.
$$
Hence we have to prove the existence of at least one
$C_{\lambda_*}(\R)$-solution  of the equation
\begin{equation}\label{bso}
z'(t) = - z(t) + q(t)z(t-h) +q(t)d(t-h).
\end{equation}
First, notice that all solutions of (\ref{bso}) are bounded on the
positive semi-axis $\R_+$ due to the boundedness of $q(t)d(t-h)$
and the exponential stability of the homogeneous $\omega$-limit
equation $z'(t) = - z(t) + g'(K)z(t-h).$  Here we use the
persistence of exponential stability under small bounded
perturbations (e.g. see \cite[Section 5.2]{CL} or \cite[Chapter VI
(9c)]{kn}) and the fact that $q(+\infty) = g'(K)$. Furthermore,
since every solution $z$ of (\ref{bso}) satisfies $z'(t) = -z(t) +
g'(K)z(t-h) + g'(K)d(+\infty) + \epsilon(t)$ with
$\epsilon(+\infty) = 0$, we get  $z(+\infty) =d(+\infty)
g'(K)(1-g'(K))^{-1}$.  Next, by effecting the change of variables
$z(t)=\exp (\lambda_* t)y(t)$ to Eq. (\ref{bso}), we get a linear
inhomogeneous equation of the form
\begin{equation}\label{bso0}
y'(t) = -(1+\lambda_*) y(t) + [p\exp(-\lambda_*h) +
\epsilon_1(t)]y(t-h) + \epsilon_2(t),
\end{equation}
where $\epsilon_1(-\infty) =0$ and $\epsilon_2(t) = O(1)$ at $t =
- \infty$. Since the $\alpha$-limit equation $y'(t) =
-(1+\lambda_*) y(t) + p\exp(-\lambda_*h)y(t-h)$ to the homogeneous
part of (\ref{bso0}) is hyperbolic, due to the above mentioned
persistence of the property of exponential dichotomy, we again
conclude that Eq. (\ref{bso0}) also has an exponential dichotomy
on $\R_-$. Thus (\ref{bso0}) has a solution $y^*$ which is bounded
on $\R_{-}$ so that  $z^*(t) = \exp (\lambda_* t)y^*(t) =
O(\exp(\lambda_* t)), \ t \to -\infty,$ is a
$C_{\lambda_*}(\R)$-solution of Eq. (\ref{bso}). Now, it is
evident that $w(t)= z^*(t)- C\psi'(t) = O(\exp(\lambda_* t))$
solves (\ref{bso}) for each $C \in \R$. In consequence, $x(t)=
d(t)+ z^*(t)- C_d\psi'(t) = ((I - \mathcal{N})^{-1}d)(t)$, if we
take $C_d = \int_{-\infty}^0(d(s)+ z^*(s))
\psi'(s))ds(\int_{-\infty}^0(\psi'(s))^2ds)^{-1}$. $\qquad \qquad
\Box$
\end{pf}
\begin{rem} \label{re1}
For $\delta >0$ small, consider $I- \mathcal{N}_1: C_{2\lambda -
\delta}(\R) \to C_{2\lambda - \delta}(\R)$, where $\mathcal{N}_1$
is defined by $ (\mathcal{N}_1x)(t) =p
\int^t_{-\infty}e^{-(t-s)}x(s-h)ds$ (recall here the discussion
after formula (\ref{llps})). Replacing $\mathcal{N}$ by
$\mathcal{N}_1$ in the proof of Lemma \ref{6}, we establish
similarly that $I- \mathcal{N}_1$ is an isomorphism of the Banach
space $C_{2\lambda - \delta}(\R)$ onto itself. Since the linear
equation $x'(t) = - x(t) + px(t-h)$ is hyperbolic, this situation
is actually simpler than the one considered in Lemma \ref{6}.
\end{rem}
\vspace{-4mm}
\section{ Existence of a continuous family of positive travelling waves}
\vspace{-3mm} In this section, we are looking for travelling waves
for (\ref{17}), that is,  solutions $u(x,t) = \phi(\varepsilon \nu
\cdot x+t), \ x, \nu \in \R^m, \ \|\nu\| = 1$, where
$c=1/\varepsilon$ is the wave speed, connecting the two equilibria
of (\ref{17}). We will suppose that $\varepsilon$ is sufficiently
small. This leads us to the question about the existence of
heteroclinic solutions to the singularly perturbed equation
\begin{equation}\label{twe}
\varepsilon^2x''(t) - x'(t)-x(t)+ g(x(t-h))=0, \quad t \in \R ,
\end{equation}
with $x(-\infty)=0,\ x(+\infty)=K$. Being a bounded function, each
travelling wave should satisfy the following integral equation
\begin{equation}\label{psea} \hspace{-7mm}
x(t) = \frac{1}{\sigma(\varepsilon)}\left(
\int_{-\infty}^te^{\frac{-2(t-s)}{1+
\sigma(\varepsilon)}}g(x(s-h))ds + \int^{+\infty}_te^{\frac{(1+
\sigma(\varepsilon))(t-s)}{2\varepsilon^2}}g(x(s-h))ds \right),
\end{equation}
where $\sigma(\varepsilon) = \sqrt{1+4\varepsilon^2}.$ For
solutions in $C_{\lambda_*}(\R)$ with $\lambda_* \in (0,\lambda)$
close to $\lambda$, this equation can be written in the shorter
form
\begin{equation}\label{be}
x - (\mathcal{I}_{\varepsilon}\circ \mathcal{G})x = 0,
\end{equation}
where $\mathcal{I}_{\varepsilon}, \mathcal{G}: C_{\lambda_*}(\R)
\to C_{\lambda_*}(\R)$ are defined by
$$
(\mathcal{I}_{\varepsilon}x)(t)=
\frac{1}{\sigma(\varepsilon)}\left(
\int_{-\infty}^te^{\frac{-2(t-s)}{1+ \sigma(\varepsilon)}}x(s-h)ds
+ \int^{+\infty}_te^{\frac{(1+
\sigma(\varepsilon))(t-s)}{2\varepsilon^2}}x(s-h)ds \right),
$$
and $ (\mathcal{G}x)(t) = g(x(t))$ is the Nemitski operator. (For
the sake of simplicity, we write $\mathcal{I}_{\varepsilon},
\mathcal{G}$ instead of $\mathcal{I}_{\varepsilon,\lambda_*},
\mathcal{G}_{\lambda_*}$). Throughout all this section, we will
suppose that the $C^1$-smooth function $g$ is defined and bounded
on the whole real axis $\R$. Clearly, this assumption does not
restrict the generality of our framework, since it suffices to
take any smooth and bounded extension on $\R_-$ of the
nonlinearity $g$ described in {\rm \bf(H)}. Notice that, since
there exists finite $g'(0)$, we have $g(x) = x\gamma(x)$ for a
bounded $\gamma \in C(\R)$. Set $ \gamma_0= \sup_{t \in
\R}|\gamma(x)|$.  As it can be easily checked, $\|\mathcal{G}x\|
\leq \gamma_0\|x\|$ so that actually $\mathcal{G}$ is
well-defined. Furthermore, we have the following lemma:
\begin{lem}\label{7} Assume that $g \in C^1(\R)$. Then $\mathcal{G}$ is Fr\'echet continuously
differentiable on $C_{\lambda_*}(\R)$ with differential
$\mathcal{G}'(x_0): y(\cdot) \to g'(x_0(\cdot))y(\cdot).$
\end{lem}
\vspace{-4mm}
\begin{pf} By the Taylor formula,
$ g(v)-g(v_0) - g'(v_0)(v-v_0) = o(v-v_0),$ $ \ v,v_0 \in \R. $
Fix some $x_0 \in C_{\lambda_*}(\R)$, then we have
$$
\|\mathcal{G}x-\mathcal{G}x_0 - g'(x_0(\cdot))(x-x_0)\| =
o(\|x-x_0\|), \quad x \in C_{\lambda_*}(\R).
$$
Clearly, it holds that $\|\mathcal{G}'(x)u\|=
\|g'(x(\cdot))u(\cdot)\| \leq \sup_{t \in \R} |g'(x(t))|\|u\|$.
Since functions in $C_{\lambda_*}(\R)$ are bounded  and $g'$ is
uniformly continuous on bounded sets of $\R$, for any given
$\delta >0$ there is $\sigma >0$ such that for $\|x-x_0\|<\sigma$
we have $ \|\mathcal{G}'(x)-\mathcal{G}'(x_0)\| <\delta. \qquad
\qquad \hspace{6.5cm} \Box $
\end{pf}
\vspace{-5mm}
 \noindent Now, we  consider the integral operators
$\mathcal{I}^+_{\varepsilon}, \mathcal{I}^-_{\varepsilon}:
C_{\mu}(\R) \to C_{\mu}(\R)$ defined as
$$
(\mathcal{I}^+_{\varepsilon}x)(t)=
 \int^{+\infty}_te^{\frac{(1+
\sigma(\varepsilon))(t-s)}{2\varepsilon^2}}x(s-h)ds,  \
(\mathcal{I}^-_{\varepsilon}x)(t)=
\int_{-\infty}^te^{\frac{-2(t-s)}{1+
\sigma(\varepsilon)}}x(s-h)ds.
$$
\begin{lem} \label{8} Set $\mathcal{I}= \mathcal{I}^-_{0}$ and $\mathcal{I}^+_{0}=0$.
If $\varepsilon \to 0+$, then $\mathcal{I}_{\varepsilon} \to
\mathcal{I}$ in the operator norm. Moreover, both operator
families $\mathcal{I}^{\pm}_{\varepsilon}: [0, 1/\sqrt{\mu}) \to
\mathcal{L}(C_{\mu}(\R),C_{\mu}(\R))$ are continuous in the
operator norm.
\end{lem}
\vspace{-5mm}
\begin{pf} We prove only that $\|\mathcal{I}_{\varepsilon}-\mathcal{I}\| \to 0$ as $\varepsilon \to 0$,
the proof of the continuous dependence of
$\mathcal{I}^{\pm}_{\varepsilon}$ on $\varepsilon$ being
completely analogous.

We first establish that $\mathcal{I}^+_{\varepsilon} \to 0$
uniformly as $\varepsilon \to 0$. In fact, for all $ t \in \R$, we
obtain
$$
|(\mathcal{I}^+_{\varepsilon}x)(t)|\leq
\int^{+\infty}_te^{\frac{t-s}{\varepsilon^2}}|x(s-h)|ds \leq
\varepsilon^2\|x\|.
$$
Furthermore, since $|x(t)|\leq \|x\|\exp(\mu t) $ for all $t \in
\R$,  for $\varepsilon^2 < 1/\mu$ we have
$$
|(\mathcal{I}^+_{\varepsilon}x)(t)|\leq
\int^{+\infty}_te^{\frac{t-s}{\varepsilon^2}}|x(s-h)|ds \leq
 \frac{\|x\|}{\varepsilon^{-2}-\mu}e^{\mu (t-h)}.
$$
Hence, for $\varepsilon^2 < (1-0.5e^{-\mu h})/\mu$, we obtain
that $\|\mathcal{I}^+_{\varepsilon}x\| \leq 2\varepsilon^2\|x\|$.

Next, we prove that $\mathcal{I}^-_{\varepsilon} \to \mathcal{I}$
uniformly as $\varepsilon \to 0$. We have
$$
|((\mathcal{I}^-_{\varepsilon} - \mathcal{I})x)(t)|\leq
\int_{-\infty}^te^{-(t-s)}(e^{\frac{\sigma(\varepsilon)-1}{
\sigma(\varepsilon)+1}(t-s)}-1)|x(s-h)|ds.
$$
Thus, for $t \leq h$, we obtain that $ \quad
|((\mathcal{I}^-_{\varepsilon} - \mathcal{I})x)(t)| \leq$
$$
\leq \int_{-\infty}^te^{-(t-s)}(e^{\frac{\sigma(\varepsilon)-1}{
\sigma(\varepsilon)+1}(t-s)}-1)\|x\|^-_{\mu}e^{\mu(s-h)}ds =
\frac{\|x\|^-_{\mu}e^{\mu(t-h)}(\sigma(\varepsilon)-1)}{
(2+(\sigma(\varepsilon)+1)\mu)(1+\mu)},
$$
and, for all $t$,
$$
|((\mathcal{I}^-_{\varepsilon} - \mathcal{I})x)(t)|\leq
\int_{-\infty}^te^{-(t-s)}(e^{\frac{\sigma(\varepsilon)-1}{
\sigma(\varepsilon)+1}(t-s)}-1)\|x\|ds =
\|x\|\frac{\sigma(\varepsilon)-1}{2}.
$$
Thus $\|\mathcal{I}^-_{\varepsilon} - \mathcal{I}\|\leq
0.5(\sigma(\varepsilon)-1)$, and the proof of the lemma is
complete. $\Box$
\end{pf}
To prove the main result of this section, stated  below as Theorem
\ref{mr}, we will make use of the following proposition:
\begin{lem}\label{lem}
Let $\{z_{\alpha}, \alpha \in A\}$, where $\N \cup \{ \infty
\}\subset A$, denote the (countable) set of roots to the equation
\begin{equation}\label{char2}
\varepsilon^2z^2-z-1+p\exp(-z h)=0.
\end{equation}
If $p> 1, \ h >0, \ \varepsilon  \in (0,1/(2\sqrt{p-1}))$ then
(\ref{char2}) has exactly two real roots $\lambda_1(\varepsilon),
\lambda_{\infty}(\varepsilon)$ such that
$$0<\lambda<  \lambda_1(\varepsilon) < 2(p-1) <
\varepsilon^{-2}-2(p-1) < \lambda_{\infty}(\varepsilon) <
\varepsilon^{-2}+1.$$ Moreover: (i) there exists an interval
$\mathcal{O} = (0,a(p,h))$ such that, for every $\varepsilon \in
\mathcal{O}$, all roots $\lambda_{\alpha}(\varepsilon), \alpha \in
A$ of (\ref{char2}) are simple  and the functions
$\lambda_{\alpha}: \mathcal{O} \to \C$ are continuous; (ii) we can
enumerate $\lambda_j(\varepsilon), j \in \N$, in such a way that
there exists $\lim_{\varepsilon \to 0+}\lambda_j(\varepsilon)
=\lambda_j$ for each $j \in \N$, where $\lambda_j \in \C$ are the
roots of (\ref{char}), with $\lambda_1=\lambda$; (iii) for all
sufficiently small $\varepsilon$, every vertical strip $\xi \leq
\Re z \leq 2(p-1)$ contains only a finite set of $m(\xi)$ roots
(if $\xi \not\in \{\Re \lambda_j, \ j \in \N\}$, then $m(\xi)$
does not depend on $\varepsilon$) $\lambda_1(\varepsilon), \dots,
\lambda_{m(\xi)}(\varepsilon)$ to (\ref{char2}), while the
half-plane $\Re z > 2(p-1)$ contains only the root $
\lambda_\infty(\varepsilon)$.
\end{lem}

\vspace{-5mm}

\begin{pf} The existence of real roots $\lambda_1(\varepsilon), \lambda_\infty(\varepsilon)$ satisfying
$\lambda< \lambda_1(\varepsilon) < \lambda_{\infty}(\varepsilon)$
is obvious when $\varepsilon  \in (0,0.5/\sqrt{p-1})$. On the
other hand, if $z_0
>0$ is a real root of (\ref{char2}), then
$ \varepsilon^2z_0^2-z_0-1<0, \ \varepsilon^2z_0^2-z_0-1+p>0. $
Hence $z_0 < (1+ \sqrt{1+4\varepsilon^2})/(2\varepsilon^2)<
\varepsilon^{-2}+1$, from which it can be checked easily that
$$
\lambda_{\infty} >
\frac{1+\sqrt{1-4(p-1)\varepsilon^2}}{2\varepsilon^2}
> \frac{1-2(p-1)\varepsilon^2}{\varepsilon^2}, \ \lambda_1
<\frac{1-\sqrt{1-4(p-1)\varepsilon^2}}{2\varepsilon^2}<2(p-1).
$$
We also notice that every multiple root $z_0$ has to satisfy the
system
\begin{equation}\label{ep} \varepsilon^2z_0^2-z_0-1+p\exp(-z_0
h)=0, \ 2\varepsilon^2z_0-1-ph\exp(-z_0 h)=0,
\end{equation}
which implies
\begin{equation}\label{ep1}
(\varepsilon^2z_0^2-z_0-1)h+2\varepsilon^2z_0-1=0, \ p\exp(-z_0
h)= \frac{2+z_0}{2+ hz_0}.
\end{equation}
The first equation of (\ref{ep1}) implies that $z_0$ is real while
the second equation of (\ref{ep}) says that $z_0 > 0$. Since $z_0$
is positive, from the first equation of (\ref{ep1}) we obtain
$0.5\varepsilon^{-2}< z_0$ (we recall that
$\varepsilon^2z_0^2-z_0-1<0$). Let $\zeta_0(p,h)$ be the maximal
positive root of the second equation of (\ref{ep1}). If
$\varepsilon
>0$ is so small that $0.5\varepsilon^{-2}> \zeta_0(p,h)$,
system (\ref{ep}) can not have any positive solution. In
consequence, the second assertion of this lemma holds if we set
$a(p,h) = 1/\sqrt{2 \zeta_0(p,h)}$.

Finally, we prove that the half-plane $\Re z
> 2(p-1)$ contains only the root $ \lambda_\infty(\varepsilon)$ of
(\ref{char2}). For this, let us evaluate $|\varepsilon^2z^2-z-1|$
on the boundary of some rectangle $[2(p-1),b]\times[-c,c] \subset
\C$, with $b,c$ being sufficiently large. For $\mu(\varepsilon),\nu(\varepsilon)$ the (real) roots of $\varepsilon^2z^2-z-1=0$, we have that
$$
|\varepsilon^2z^2-z-1| =
\varepsilon^2|z-\mu(\varepsilon)||z-\nu(\varepsilon)| \geq
\varepsilon^2|\Re z-\mu(\varepsilon)||\Re z-\nu(\varepsilon)| =
|\varepsilon^2(\Re z)^2-\Re z-1|.
$$
Thus, for $\Re z = 2(p-1)$, we obtain
$$|\varepsilon^2 z^2-z-1|  \geq \Re z+1 - \varepsilon^2(\Re
z)^2
> p.
$$
If $\Re z > 2(\varepsilon^{-2}+1)$, then
$$|\varepsilon^2z^2-z-1|  \geq \varepsilon^2(\Re z)^2 -\Re
z-1
> 8p -3 > p.
$$
Similarly, for $|\Im z| > p/\varepsilon$ fixed, we get
$$
|\varepsilon^2z^2-z-1| =
\varepsilon^2|z-\mu(\varepsilon)||z-\nu(\varepsilon)| \geq
\varepsilon^2(\Im z)^2 > p.
$$
Thus, by Rouch\'e's theorem, $\varepsilon^2z^2-z-1+p\exp(-z h)=0$
and  $\varepsilon^2z^2-z-1=0$ have the same number of roots in the
half-plane $\Re z > 2(p-1)$, that is exactly one root.

Therefore, for all $\lambda_j$ with $\Re \lambda_j \in [\xi,
2(p-1)]$ and $\varepsilon  \in (0,0.25/\sqrt{p-1})$, we get
$$
p e^{-\xi h} \geq |\Im(\varepsilon^2\lambda_j^2-\lambda_j-1)| =
|\Im \lambda_j||1-2\varepsilon^2\Re \lambda_j| \geq |\Im
\lambda_j|/2,
$$
so that $|\Im \lambda_j| \leq 2pe^{-\xi h}$. Hence, applying
Rouch\'e's theorem to the functions $\varepsilon^2z^2-z-1+p\exp(-z
h)$ and $-z-1+p\exp(-z h)$ along an appropriate rectangle inside
 $[\xi-1, 2(p-1)]\times [-3pe^{-\xi h}, 3pe^{-\xi h}] \subset \C$, we
prove the last assertion of Lemma \ref{lem}. $\qquad \Box$
\end{pf}
\begin{thm} \label{mr}  Assume {\rm \bf(H)} and that the positive equilibrium of
Eq. (\ref{17a}) is globally attractive. Let  $\psi $ be some
heteroclinic orbit of Eq. (\ref{17a}): $\psi(-\infty) = 0,$ $
\psi(+\infty) = K$.  Then, for every $\delta >0$ there is an
interval $\mathcal{E}= (-\varepsilon_0,\varepsilon_0)$ and a
continuous family of positive heteroclinic orbits
$\psi_\varepsilon: \mathcal{E} \to C_{\lambda - \delta}(\R)$ of
Eq. (\ref{twe}) satisfying the additional conditions $\psi_0 =
\psi$ and $\int_{-\infty}^0\psi_\varepsilon(s)\psi'(s)ds =
0.5\psi^2(0)$. Furthermore, for every $\varepsilon \in
\mathcal{E}\setminus\{0\}$ we have that $\psi_\varepsilon(t-t_0) =
\exp (\lambda_1(\varepsilon) t) + O(\exp(2\lambda t))$ at $t \to -
\infty$ for some $t_0 = t_0(\varepsilon)\in \R$, and that
$\psi'_\varepsilon(t-t_0) = \lambda_1(\varepsilon) \exp
(\lambda_1(\varepsilon) t) + O(\exp(2\lambda  t))> 0$ on some
semi-axis $(-\infty, z]$.
\end{thm}
\begin{pf} We represent the mentioned orbit
$\psi$ of (\ref{17a}) as $\psi  =  \alpha \psi' + \phi_0$, where
$$\phi_0 = (\psi - \alpha \psi')\in C_{\psi,\lambda_*}(\R),  \quad
\alpha = \psi^2(0)(2\int_{-\infty}^0(\psi'(s))^2ds)^{-1}.
$$
For $\delta >0$ small, consider $\lambda_*=\lambda -\delta$.
In virtue of Lemmas \ref{6}, \ref{7} and \ref{8}, we can apply the
implicit function theorem (e.g., see \cite[pp. 36-37]{ap} or
\cite[p. 170]{smo}) to Eq. (\ref{be}) written as $F(\phi,
\varepsilon) =0$, where $F: C_{\psi,\lambda_*}(\R)\times \R \to
C_{\lambda_*}(\R)$,
$$
F(\phi, \varepsilon) =\alpha \psi' + \phi -
(\mathcal{I}_{\varepsilon}\circ \mathcal{G})(\alpha \psi' + \phi),
\ {\rm and \ } \mathcal{I}_{0} = \mathcal{I}.
$$
Observe that $F(\phi_0, 0) =0$ and $F_{\phi}(\phi_0, 0) = I -
\mathcal{N}$. In this way, we establish the existence of an
interval $\mathcal{E}= (-\varepsilon_0,\varepsilon_0), \
\varepsilon_0  \in (0,1/(2\sqrt{p-1}))$ and a continuous family
$\phi_{\varepsilon} : \mathcal{E}\to C_{\psi,\lambda_*}(\R)$ of
solutions to $F(\phi, \varepsilon) =0$. Notice that $\psi_0=\psi$,
$\psi_\varepsilon=\alpha \psi' + \phi_\varepsilon \in
C_{\lambda_*}(R)$ satisfy Eq. (\ref{be}), so that, as it can be
checked directly, $\psi_{\varepsilon}(+\infty) =
g(\psi_{\varepsilon}(+\infty))$. Thus $\psi_{\varepsilon}(+\infty)
=K$ and $\psi_\varepsilon$ satisfies all conclusions of the third
sentence of the theorem, except its positivity,
 which is proved  below.

Assume now that $\varepsilon_0$ is sufficiently small so that
$\lambda_1(\varepsilon) < 0.5\lambda_\infty(\varepsilon)$ for all
$\varepsilon \in \mathcal{E}\setminus \{ 0\}$. Since $g(x) = px + O(x^2)$ as $x
\to 0$, and since there exists a constant $C_1 >0$ such that
$|\psi_{\varepsilon}(t)| \leq C_1 \exp(\lambda_* t), \ t \leq 0,\
\varepsilon \in \mathcal{E}$,
we get
\begin{equation}\label{pse}
\varepsilon^2 \psi_{\varepsilon}''(t)- \psi_{\varepsilon}'(t)  -
\psi_{\varepsilon}(t) + p \psi_{\varepsilon} (t-h) =
\Psi_{\varepsilon}(t),
\end{equation}
where  $\Psi_{\varepsilon}(\cdot)=p  \psi_{\varepsilon}
(\cdot-h)-g(\psi_{\varepsilon} (\cdot-h)) \in C_{2\lambda_*}(\R)$.
Moreover, $\|\Psi_{\varepsilon}\|_{2\lambda_*} \leq C_2$ for some
$C_2 >0$ which does not depend on $\varepsilon$. Now,
$C_{2\lambda_*}(\R)$-solutions $x_{\varepsilon}$ to
\begin{equation}\label{lpse}
\varepsilon^2 x''(t) - x'(t)  - x(t) + px(t-h) =
\Psi_{\varepsilon}(t),
\end{equation}
are solutions to the equation $(I-
p\mathcal{I}_{\varepsilon})x_{\varepsilon} =
-\mathcal{I}_{\varepsilon}\Psi_{\varepsilon}$. Due to Remark
\ref{re1} and Lemma \ref{8}, for $\lambda_* = \lambda
- \delta$ close to $\lambda$  the operator $I-
p\mathcal{I}_{\varepsilon}$ is invertible in $C_{2\lambda_*}(\R)$
for all sufficiently small $\varepsilon$. Moreover, Lemma \ref{8}
implies that there exists a subinterval $\mathcal{E}_1 \subset
\mathcal{E}$ such that $\|(I- p\mathcal{I}_{\varepsilon})^{-1}\|
\leq C_3$ for all $\varepsilon \in \mathcal{E}_1$. Hence, we
obtain $\|x_{\varepsilon}\| \leq \|(I-
p\mathcal{I}_{\varepsilon})^{-1}\mathcal{I}_{\varepsilon}\|\|\Psi_{\varepsilon}\|
\leq C_4$ for all $\varepsilon \in \mathcal{E}_1$. Therefore Eq.
(\ref{lpse}) has a bounded solution $x_{\varepsilon}$ such that
$|x_{\varepsilon}(t)| \leq C_4 \exp(2\lambda_* t), \ t \leq 0,\
\varepsilon \in \mathcal{E}_1$. Consequently, $z_{\varepsilon}(t)
= \psi_\varepsilon(t) - x_{\varepsilon}(t)$ solves the linear
homogenous equation
$$
\varepsilon^2 z''(t) - z'(t)  - z(t) + pz(t-h) = 0, \quad t \in \R
$$
and is bounded as $t \to - \infty$. This is possible if and only
if
$$
z_{\varepsilon}(t) = A_{\varepsilon} \exp
(\lambda_1(\varepsilon)t) + B_{\varepsilon} \exp
(\lambda_\infty(\varepsilon)t) + \sum_{j=2}^NC_{j,\varepsilon}
\exp (\lambda_j(\varepsilon) t),
$$
where $\lambda_j(\varepsilon) \in \C, \ j = 1, \dots, N$ and
$\lambda_\infty(\varepsilon)$ are the roots with non-negative real
parts of the characteristic equation (\ref{char2}), $A_\varepsilon , B_\varepsilon \in\R, C_{j,\varepsilon}\in\C$, $\varepsilon\in\mathcal{E}_1\setminus \{ 0\}$. In
consequence,
$$
\psi_\varepsilon(t) = A_{\varepsilon} \exp
(\lambda_1(\varepsilon)t) + B_{\varepsilon} \exp
(\lambda_\infty(\varepsilon)t) + \sum_{j=2}^NC_{j,\varepsilon}
\exp (\lambda_j(\varepsilon) t)+ x_{\varepsilon}(t).
$$
It follows from Lemma \ref{lem} that $\Re \lambda_j(\varepsilon) <
\lambda_* < \lambda < \lambda_1(\varepsilon)\ <
0.5\lambda_\infty(\varepsilon)$, provided $\varepsilon$ is small
(say, $\varepsilon \in \mathcal{E}_2 \subset \mathcal{E}_1$) and
$\lambda_*$ is sufficiently close to $\lambda$.  Since
$\psi_\varepsilon (t) = O(\exp(\lambda_* t))$, this implies
immediately that  $C_{j,\varepsilon}=0$  and
$$
\psi_\varepsilon (t) =  A_{\varepsilon} \exp
(\lambda_1(\varepsilon)t) + B_{\varepsilon} \exp
(\lambda_\infty(\varepsilon)t) + x_{\varepsilon}(t), \ t \in \R, \
\varepsilon \in \mathcal{E}_2\setminus \{ 0\}.
$$
To prove the positivity of $\psi_\varepsilon$ for $\varepsilon$
small, we first establish that $\limsup_{\varepsilon \to
0}|B_\varepsilon|$ is finite, from which we  deduce that the
constants $A_\varepsilon$ are positive; in fact, we will find that
$A_\varepsilon > 1 - 4\delta$. Let us suppose already that
$1-5 \delta >0$ and $\lambda_1(\varepsilon)<2\lambda_*
-\delta$ for all $\varepsilon \in \mathcal{E}_2$.
 Since $\psi_\varepsilon \in C_{\lambda_*}(R)$, for all $t\le 0,\
\varepsilon \in \mathcal{E}_2\setminus \{ 0\}$, we get
$$|A_\varepsilon \exp(\lambda_1(\varepsilon )t)+
B_\varepsilon \exp(\lambda_\infty (\varepsilon )t)| \le
|\psi_\varepsilon (t)|+|x_\varepsilon (t)|\le C_5\exp(\lambda_*t),
$$
where $C_5=C_1+C_4$. In particular, taking $t=0$ and $t=-1$, we
obtain
$$|A_\varepsilon +B_\varepsilon|\leq C_5,
|A_\varepsilon +B_\varepsilon \exp(\lambda_1(\varepsilon
)-\lambda_\infty (\varepsilon )) |\leq C_5\exp(
\lambda_1(\varepsilon )-\lambda_*)\leq C_5\exp(\lambda),$$ hence
$|B_\varepsilon |(1-\exp(\lambda_1(\varepsilon )-\lambda_\infty
(\varepsilon )))\le C_5(1+\exp(\lambda)) := C_6,$ for  $\varepsilon
\in \mathcal{E}_2\setminus \{ 0\}$.

Noting that $\exp(\lambda_1(\varepsilon )-\lambda_\infty
(\varepsilon ))\to 0$ as $\varepsilon \to 0$, we deduce that there
is $\mathcal{E}_3 = (-\varepsilon_3,\varepsilon_3)\subset
\mathcal{E}_2$   such that $|B_\varepsilon |\le 2C_6$ for
$\varepsilon \in \mathcal{E}_3\setminus \{ 0\}$, so that
$$|B_\varepsilon \exp(\lambda_\infty (\varepsilon )t)|\leq 2C_6 \exp(\lambda_\infty
(\varepsilon )t) \le 2C_6 \exp(\delta t)\exp((2\lambda_* -\delta
)t),$$ for $ t\leq 0, \ \varepsilon \in \mathcal{E}_3\setminus \{
0\}.$ Set $y_\varepsilon (t)=B_\varepsilon \exp(\lambda_\infty
(\varepsilon )t)+x_\varepsilon (t)$. By Lemma \ref{psi}, we have
$\psi (t)=\exp(\lambda t)+z(t)$ with $z(t)=O(\exp(2\lambda_* t))$
at $t=-\infty$. Since $\lim_{t \to -\infty}C_6 \exp(\delta t) =0$,
we now conclude that there is $s_0=s_0(\delta)< 0$ such that for
$t\leq s_0$ and $0<|\varepsilon |< \varepsilon_3$ we have
$$|y_\varepsilon (t)|\leq \delta \exp((2\lambda_* -\delta)t), \
|y_\varepsilon (t)-z(t)|\leq \delta \exp((2\lambda_* -\delta)t).
$$ On the other hand, for  $\delta_0=\delta \exp((\lambda
-\lambda_*)s_0)=\delta \exp(\delta s_0)$, there exists
$\varepsilon_4=\varepsilon_4(\delta )\in (0,\varepsilon_3]$ such
that, for $|\varepsilon | < \varepsilon_4$, we have
$|\psi_\varepsilon (t)-\psi (t)|\leq \delta_0 \exp(\lambda_* t). $
Taking  $t=s_0$   we obtain for $0<|\varepsilon | < \varepsilon_4$
$$
|A_\varepsilon \exp(\lambda_1(\varepsilon ) s_0)-\exp(\lambda
s_0)|\leq |\psi_\varepsilon (s_0)-\psi (s_0)|+|y_\varepsilon
(s_0)-z(s_0)| \leq 2\delta\ \exp(\lambda s_0),$$ hence
$\psi_\varepsilon (s_0)=A_\varepsilon \exp(\lambda_1(\varepsilon)
s_0)+y_\varepsilon (s_0)\geq \exp(\lambda s_0)-(|A_\varepsilon
\exp(\lambda_1(\varepsilon ) s_0)-\exp(\lambda
s_0)|+|y_\varepsilon (s_0)|)>(1-3\delta )\exp(\lambda s_0)$.
Therefore, for all $0< |\varepsilon |< \varepsilon_4$,
$$
 A_\varepsilon > -y_\varepsilon (s_0)\exp(-\lambda_1(\varepsilon )s_0)+(1-3\delta)
\exp((\lambda -\lambda_1(\varepsilon ))s_0)
 \geq (1-4\delta )>0.
 $$
Thus, for $0<|\varepsilon|< \varepsilon_4$ and $t\leq s_0$ we get
$\psi_\varepsilon (t)\geq \exp(\lambda_1(\varepsilon )t)
[(1-4\delta)-\delta )]>0.$ Since $\lim_{\varepsilon \to
0}\psi_\varepsilon =\psi $ uniformly on $\R$ and $\psi$ is bounded
from below by a positive constant on  $[s_0,\infty)$, we conclude
that $\psi_\varepsilon$ is positive on $\R$, for all $\varepsilon
$ small.

Finally, for every fixed $\varepsilon \in \mathcal{E}_2 \setminus \{ 0\}$, we have
that $ g(x) = px + q(x)x^2,$ $ \psi_\varepsilon(t) =
A_{\varepsilon} \exp (\lambda_1(\varepsilon)t) +
b_{\varepsilon}(t)\exp(2\lambda_*t),$ where $q \in C[0,+\infty)$
and $b_{\varepsilon}$ is bounded on $(-\infty,0]$. Hence,
$\psi_\varepsilon(t) = A_\varepsilon \exp
(\lambda_1(\varepsilon)t) + O(\exp(2\lambda_* t))$ at $-\infty$
and
$$
g(\psi_\varepsilon(t-h)) = A_\varepsilon
p\exp(-\lambda_1(\varepsilon)h) \exp (\lambda_1(\varepsilon)t) +
c_\varepsilon(t)\exp(2\lambda_* t),\ \varepsilon \in
\mathcal{E}_2\setminus\{0\},
$$
where $c_\varepsilon(t)$ is bounded: $|c_\varepsilon(t)| \leq
c_0(\varepsilon), \ t \leq 0$.
 Differentiating (\ref{psea}), we
obtain
\begin{eqnarray} \label{2f}
&  & \psi_\varepsilon'(t) =
\frac{1}{\sigma(\varepsilon)}\left(\frac{-2}{1+\sigma(\varepsilon)}
\int_{-\infty}^te^{\frac{-2(t-s)}{1+ \sigma(\varepsilon)}}g(\psi_\varepsilon(s-h))ds + \right.\\
&  & \left. +\frac{1+\sigma(\varepsilon)}{2\varepsilon^2}
\int^{+\infty}_te^{\frac{(1+\sigma(\varepsilon))(t-s)}{2\varepsilon^2}}g(\psi_\varepsilon(s-h))ds
\right)=  \nonumber \\
&  & = \frac{A_\varepsilon
p\exp(-\lambda_1(\varepsilon)h)}{\sigma(\varepsilon)}\left(\frac{-2}{1+
\sigma(\varepsilon)} \int_{-\infty}^te^{\frac{-2(t-s)}{1+
\sigma(\varepsilon)}} \exp
(\lambda_1(\varepsilon)s)ds +\right. \nonumber\\
&  & + \left. \frac{1+
\sigma(\varepsilon)}{2\varepsilon^2}\int^{+\infty}_te^{\frac{(1+
\sigma(\varepsilon))(t-s)}{2\varepsilon^2}} \exp
(\lambda_1(\varepsilon)s)ds \right) +
O(\exp(2\lambda_* t))= \nonumber\\
&   &   =  A_\varepsilon\lambda_1(\varepsilon) \exp
(\lambda_1(\varepsilon)t) + O(\exp(2\lambda_* t)),  \quad t \to -
\infty. \nonumber
\end{eqnarray}
Hence, $\psi_{\varepsilon}, \psi'_{\varepsilon} \in
C_{\lambda_1(\varepsilon)}(\R)$ so that $\Psi_{\varepsilon} \in
C_{2\lambda_1(\varepsilon)}(\R)$ in (\ref{pse}).
 Therefore, in view of \cite[Proposition 7.1]{FA} and
the inequality $\lambda< \lambda_1(\varepsilon), \ \varepsilon \in
\mathcal{E}_2\setminus\{0\}$, we get from (\ref{pse}) that
$\psi_\varepsilon(t) = A_{\varepsilon} \exp
(\lambda_1(\varepsilon)t) + O(\exp(2\lambda t))$ and
$\psi'_\varepsilon(t) = A_\varepsilon\lambda_1(\varepsilon) \exp
(\lambda_1(\varepsilon)t) + O(\exp(2\lambda t)), \quad t \to
-\infty$. $\qquad \qquad  \qquad \qquad \qquad \qquad \Box $
\end{pf}
\section{Non-monotonicity of travelling waves}
As it was noticed in \cite{gous,gou,gouk,gouss,lw,swz}, various
investigators have studied numerically the case of a large delay
in the Nicholson's blowflies equation, and noted a loss of
monotonicity of the wave front as the delay  increases, {\it
``with the front developing a prominent hump"} whose height {\it
``is bounded above by a bound that does not depend on the delay"},
see \cite[p. 308]{gous} from which the above citation was taken.
It is not difficult to explain the second phenomenon, since at
every point of local maximum  $\sigma$ of $\psi(t,c)$ we have
$\psi'(\sigma,c) = 0, \ \psi''(\sigma,c) \leq 0$ so that $
\psi(\sigma,c) \leq g(\psi(\sigma-h,c)) \leq \max_{x \geq 0}
g(x)$. Here we explain also the first phenomenon, easily getting
the oscillation of the travelling waves  about $K$ as $t\to +\infty$ stated in  Theorem \ref{main} from the next two
lemmas.
\begin{lem} \label{15}
Let $g'(K) <0, h > 0$ and $|g'(K)|he^{h+1} >1$.  Then the
 equation
\begin{equation}\label{char5}
\varepsilon^2z^2-z-1+g'(K)\exp(-z h)=0
\end{equation}
has no negative real roots, for all sufficiently small
$\varepsilon$. Moreover, if the equilibrium $K$ of (4) is
hyperbolic, then, for all small $\varepsilon$, there are no roots
of (28) on the imaginary axis.
\end{lem}
\begin{pf} Set $\Delta_\varepsilon (z)=\varepsilon^2z^2-z-1+g'(K)\exp(-z h)$.
We first prove that the lemma is valid for $\varepsilon =0$ (see
also \cite{gl}). Let $z_0$ be the maximum point of $\Delta_0(z)$
on $\R$, i.e., $z_0\in \R$ is such that
$\Delta_0'(z_0)=-1-hg'(K)\exp(-z_0 h)=0$. Note that
$\Delta_0(-\infty)=-\infty$ and $\Delta_0(z)<0$ for $z\geq 0$. If
there is a negative zero of $\Delta_0(z)$, then $z_0<0$ and
$\Delta_0(z_0)=-z_0-1-1/h\geq 0$, implying that
$0=\Delta_0'(z_0)\geq -1+|g'(K)|he^{h+1}$, which contradicts the
hypothesis $|g'(K)|he^{h+1} >1$.

If $h|g'(K)|\geq 1$, then $\Delta_\varepsilon '(z)\geq
2z\varepsilon^2-1+e^{-zh}>z(2\varepsilon^2-h)>0$ for all $z<0$ and
$\varepsilon^2<h/2$, hence $\Delta_\varepsilon (z)<0$ for $z\leq
0$. Now, let $h|g'(K)|<1$, so that $z_0<0$. For $|\varepsilon|>0$
small, by the implicit function theorem we conclude that there is
a  negative root $z(\varepsilon)$ of the equation
$\Delta_\varepsilon'(z)=0$ with $z(0)=z_0$; moreover,
$z(\varepsilon)$ is the absolute maximum point of $z\mapsto
\Delta_\varepsilon(z)$ on $(-\infty, 0]$. Since
$\delta(\varepsilon):=\Delta_\varepsilon (z(\varepsilon))$ depends
continuously on $\varepsilon$ and $\delta (0)<0$,
 for $\varepsilon >0$ small we have
$\Delta_\varepsilon(z)<0$ for all $z\leq0$.

We now prove that (\ref{char5}) has no roots on the imaginary
axis. First, notice that $|\Delta_\varepsilon (ib)| \geq |g'(K)|
>0$ if $b >2|g'(K)|$.  For
$\varepsilon =0$, Eq. (\ref{char5}) does not have roots on the
imaginary axis, therefore $|\Delta_0(ib)|>0, \ |b| \leq 2|g'(K)|$.
Hence, $|\Delta_\varepsilon (ib)|>0$ for all $\varepsilon $ small
and $|b| \leq 2|g'(K)|$, which implies the hyperbolicity of Eq.
(\ref{char5}). $\Box$
\end{pf}

The next lemma can be considered as an extension of the linearized
oscillation theorem from \cite{gt} to the second order delay
differential equation
\begin{equation}\label{twee}
\varepsilon^2x''(t) - x'(t)-x(t)+ g(x(t-h))=0, \quad t \in \R.
\end{equation}
\begin{lem} \label{lll} Assume {\rm
\bf(H)} and that $g'(K)he^{h+1} <-1$. For  small $\varepsilon >0,$
set $(Dx)(t) = \varepsilon^2x''(t) - x'(t)-x(t)+ g'(K)x(t-h)$.
Then every non constant and bounded solution $x:\R \to \R$ of
$(\ref{twee})$ such that $x(+\infty) =K$ oscillates about $K$.
\end{lem}
\begin{pf} Consider some non-constant solution $x:\R \to \R$ of
(\ref{twee}) such that $x( +\infty) =K$. If we suppose for a
moment that, for some $\eta \in \R$, it holds $x(s) = K$
identically for all $s \geq \eta$,  then we obtain easily that
$x(s) = K$ for all $s \in [\eta - h,\eta]$. Hence $x$ should be a
constant solution, in contradiction with our initial assumption.
Therefore either  $\sigma(t) = x(t)-K$ oscillates about zero or is
eventually non-constant and non-negative, or non-positive. In order to get a contradiction, assume that $\sigma$ is
not oscillatory. Notice that $\sigma$ satisfies the following
linear asymptotically autonomous delay differential equation
\begin{equation}\label{tweee}
\hspace{-7mm} \varepsilon^2\sigma''(t) - \sigma'(t)-\sigma(t)+
\gamma(t)\sigma(t-h)=0,\  t \in R,\ \gamma(+\infty) = g'(K) < 0,
\end{equation}
where $\gamma(t) = g'(K) + c_0(t)$ and
$c_0(t)=g'(K+\theta(t)\sigma(t-h))-g'(K)$ for some $\theta(t)\in
(0,1)$ given by the mean value theorem. Since $x(t)$ is bounded on
$\R$ and $x(t) \to K$ as $t \to +\infty$, we can use the integral
representations (\ref{psea}) and (\ref{2f}) to prove that $\lim_{t
\to +\infty} x'(t)=0$. From Lemma \ref{15}, it follows that  the
equation $(Dx)(t)=0$ is hyperbolic, hence  the equilibrium $(K,0)$
of the system $x'(t)=v(t),\varepsilon^2 v'(t)
-v(t)-x(t)+g(x(t-h))=0$ is hyperbolic for all sufficiently small
$\varepsilon$. Thus the trajectory of $x(t)$ belongs to the stable
manifold of the hyperbolic equilibrium $K$ of (\ref{twee}), so
that we can find $a
>0$ such that $\sigma(t) =O(e^{-at}), \sigma'(t) =O(e^{-at})$
at $t=+\infty$. Therefore we have
$c_0(t)=O(\sigma(t-h))=O(e^{-at})$ at $+\infty$.  From \cite[Proposition 7.2]{FA}
 (see also \cite[Proposition 2.2]{hl1}), we
conclude that: $(i)$ either there are $b\geq a,\ \delta>0$ and
$u(t)$ a nontrivial eigensolution of the limiting equation
\begin{equation}\label{30}
\varepsilon^2u''(t)-u'(t)-u(t)+g'(K)u(t-h)=0
\end{equation}
on the generalized eigenspace associated with the (nonempty) set
$\Lambda$ of eigenvalues  with $\Re e\, \lambda=-b$, such that
$\sigma(t)=u(t)+O(\exp({-(b+\delta)t}),\quad t\to+\infty;$ \\
$(ii)$ or $\sigma(t)$ decays superexponentially at $+\infty$.
However, this latter condition is not possible: as it was
established in \cite[Lemma 3.1.1]{hl2}, if $\gamma(+\infty) \not=
0$ then every eventually nontrivial and nonnegative solution of
(\ref{tweee}) does not decay superexponentially (see also
\cite[Lemma A.1]{hl1} for the case $\gamma(+\infty) > 0$). On the
other hand,  from Lemma 15 we know that  there are no real
negative eigenvalues of (\ref{30}): hence $\Im m\, \lambda \neq 0$
for all $\lambda\in\Lambda$. From \cite[Lemma 2.3]{hl1}, we
conclude that $\sigma (t)$ is oscillatory. \qquad \qquad \qquad
\qquad \qquad \qquad \qquad \qquad \qquad \qquad \qquad$\Box$
\end{pf}
Finally, we observe that due to the exponential stability of the
positive steady state, which implies fast convergence, numerical
heteroclinic solutions $\psi(t,c)$  exhibit only one or two well
pronounced humps, see \cite[Fig. 2]{gou}.

\vspace{-5mm}

\section*{Acknowledgements} This research was supported by FONDECYT (Chile),
projects 7040044 (Teresa Faria) and 1030992 (Sergei Trofimchuk),
by FCT (Portugal), program POCTI/ FEDER, under CMAF (Teresa Faria), and
by CONICYT (Chile) through PBCT program ACT-05 (Sergei
Trofimchuk). This work was partially written while T. Faria was
visiting the University of Talca, and she thanks the University
for its kind hospitality. \vspace{-3mm}

 The authors express their gratitude to the anonymous referee, whose valuable comments helped to improve the original version of this paper.

\end{document}